\documentclass[11pt]{amsart}
\usepackage{amsxtra}
\usepackage{amssymb}
\theoremstyle{plain}
\newcommand{\cleqn}{\setcounter{equation}{0}}
\newcommand{\clth}{\setcounter{theorem}{0}}
\newcommand {\sectionnew}[1]{\section{#1}\cleqn\clth}
\newtheorem{theorem}{Theorem}[section]
\newtheorem{lemma}[theorem]{Lemma}
\newtheorem{definition-theorem}[theorem]{Definition-Theorem}
\newtheorem{proposition}[theorem]{Proposition}
\newtheorem{corollary}[theorem]{Corollary}
\newtheorem{definition}[theorem]{Definition}
\newtheorem{example}[theorem]{Example}
\newtheorem{remark}[theorem]{Remark}
\newtheorem{conjecture}[theorem]{Conjecture}
\newtheorem{notation}[theorem]{Notation}
\newcommand \bth[1] { \begin{theorem}\label{t#1} }
\newcommand \ble[1] { \begin{lemma}\label{l#1} }

\newcommand \bpr[1] { \begin{proposition}\label{p#1} }
\newcommand \bco[1] { \begin{corollary}\label{c#1} }
\newcommand \bde[1] { \begin{definition}\label{d#1}\rm }
\newcommand \bex[1] { \begin{example}\label{e#1}\rm }
\newcommand \bre[1] { \begin{remark}\label{r#1}\rm }
\newcommand \bcj[1] { \begin{conjecture}\label{j#1}\rm }

\newcommand \bnota[1] { \begin{notation}\label{n#1}\rm }
\renewcommand {\eth} { \end{theorem} }
\newcommand {\ele} { \end{lemma} }

\newcommand {\epr} { \end{proposition} }
\newcommand {\eco} { \end{corollary} }
\newcommand {\ede} { \end{definition} }
\newcommand {\eex} { \end{example} }
\newcommand {\ere} { \end{remark} }
\newcommand {\ecj} { \end{conjecture} }

\newcommand {\enota} { \end{notation} }
\newcommand \thref[1]{Theorem \ref{t#1}}

\newcommand \prref[1]{Proposition \ref{p#1}}

\newcommand \cjref[1]{Conjecture \ref{j#1}}

\newcommand \lb[1]{\label{#1}}
\def \KK {{\mathbb K}}
\def \Zset {{\mathbb Z}}

\def \ZNp  {{\mathbb Z}^N_+}
\def \Nset   {{\mathbb N}}

\def \UU {{\mathcal{U}}}

\def \HH {{\mathcal{H}}}
\def \TT {{\mathcal{T}}}

\def \qb {{\bf{q}}}
\def \db {{\bf{d}}}

\def \De {\Delta}   % Greek letters

\def \sig {\sigma}

\def \sig{\sigma}
\def \mt  {\mapsto}
\def \hra {\hookrightarrow}

\def \rcor {\rangle}
\def \lcor {\langle}
\def \ol {\overline}

\def \wh {\widehat}

\def \id { {\mathrm{id}} }

\def \g  {\mathfrak{g}}   % Lie algebra letters

\DeclareMathOperator \Aut { {\mathrm{Aut}} }

\DeclareMathOperator \Ker { {\mathrm{Ker}} }

\DeclareMathOperator \Fract { {\mathrm{Fract}} }

\renewcommand \max { {\mathrm{max}} }

\begin{document}
%%%%%%%%%%%%%%%%%%%%%%%%%%%%%%%%%%%%%%%%%%%%%%%%%%%%%%%%%%%%%%%%%%%%%%%%%%%
%%%%%%%%%%%%%%%%%%%%%%    Title    %%%%%%%%%%%%%%%%%%%%%%%%%%%%%%%%%%%%%%%%
\title[The Launois--Lenagan conjecture]
{The Launois--Lenagan conjecture}
\author[Milen Yakimov]{Milen Yakimov}
\thanks{The author was supported in part
by NSF grant DMS-1001632.}
\address{
Department of Mathematics \\
Louisiana State University \\
Baton Rouge, LA 70803
U.S.A.
}
\email{yakimov@math.lsu.edu}
\date{}
\begin{abstract} In this note we prove the Launois--Lenagan conjecture 
on the classification of the automorphism groups of the algebras 
of quantum matrices $R_q[M_n]$ of square shape
for all positive integers $n$, base fields $\KK$, and deformation
parameters $q \in \KK^*$ which are not roots of unity.
\end{abstract}
\maketitle
%%%%%%%%%%%%%%%%%%%%   Introduction   %%%%%%%%%%%%%%%%%%%%%%%%%%%%%%%%%%%%%%%%
\sectionnew{Introduction}
\lb{intro}
%%%%%%%%%%%%%%%%%%%%%%%%%%%%%%%%%%%%%%%%%%%%%%%%%%%%%%%%%%%%%%%%%%%%%%%%%%%%%
The explicit description of automorphism groups of noncommutative algebras
is a difficult problem even for algebras of low Gelfand--Kirillov dimension. Such was 
achieved in small number of examples. Andruskiewitsch--Dumas 
\cite{AD} and Launois--Lenagan \cite{LaL1} made two conjectures which 
predicted the explicit structure of the automorphism groups of important 
families of quantized universal enveloping algebras of nilpotent Lie algebras.
The algebras in both families have arbitrarily large Gelfand--Kirillov 
dimensions. 
The former conjecture concerns the positive parts of the quantized
universal enveloping algebras of all simple Lie algebras 
$\UU_q^+(\g)$ and the latter one the algebras of quantum matrices of 
square shape $R_q[M_n]$. We proved the former conjecture in \cite{Y-ad}
by exhibiting a relationship to a certain type of ``bi-finite unipotent'' 
automorphisms of completed quantum tori and proving a rigidity for those.
Here we use this method and results of Launois and Lenagan 
to settle the latter conjecture.

Let $\KK$ be an arbitrary base field, $q \in \KK^*$ an element which 
is not a root of unity, and 
$n$ a positive integer. The algebra of quantum matrices $R_q[M_n]$ is the 
$\KK$-algebra with generators $x_{kl}$, $1 \leq k, l \leq n$
and relations
\begin{align*}
x_{ij} x_{kj} &= q x_{kj} x_{ij}, \quad \mbox{for} \; i < k, \\
x_{ij} x_{il} &= q x_{il} x_{ij}, \quad \mbox{for} \; j < l, \\ 
x_{ij} x_{kl} &= x_{kl} x_{ij}, \quad 
\mbox{for} \; i < k, j > l,\\
x_{ij} x_{kl} - x_{kl} x_{ij} &= (q-q^{-1}) x_{il} x_{kj}, \quad 
\mbox{for} \; i < k, j<l.
\end{align*}
The torus
\[
\HH := (\KK^*)^{2n}/\{ (c, \ldots, c) \mid c \in \KK^*\} 
\cong (\KK^*)^{ 2n -1}
\]
acts on $R_q[M_n]$ by algebra automorphisms by 
\begin{equation}
\label{torus2}
(c_1, \ldots, c_{2n}) \cdot x_{kl} 
:= c_k c_l^{-1} x_{kl}, 
\quad (c_1, \ldots, c_{2n}) \in (\KK^*)^{2 n}.
\end{equation}
One also has the transpose automorphism 
\[
\tau \in \Aut R_q[M_n] \; \; 
\mbox{given by} \; \; 
\tau(x_{kl}) = x_{lk}.
\]
\bcj{LL-cnj} (Launois--Lenagan, \cite{LaL1})
For all base fields $\KK$, $q \in \KK^*$ 
not a root of unity, and integers $n >1$ the group of algebra 
automorphisms of $R_q[M_n]$ is isomorphic 
to $\HH \rtimes \Zset_2$ via the identification 
$\Zset_2 \cong \lcor \tau \rcor$. 
\ecj
The case of $n=2$ was proved by Alev and Chamarie \cite{ACh}
before this conjecture was stated. Recently,
Launois and Lenagan proved it for $n=3$ 
in \cite{LaL2}. Earlier they also classified in \cite{LaL1} the 
automorphism groups of quantum matrices of size $m \times n$ 
for $m \neq n$. The methods of \cite{LaL1} rely on
specific properties of the height one prime ideals in that case
and cannot be used to prove \cjref{LL-cnj}. The conjecture was 
open for all $n>3$.

The algebra $R_q[M_n]$ is a connected $\Nset$-graded algebra
generated in degree one by assigning $\deg x_{kl} = 1$, 
$\forall 1 \leq k, l \leq n$. For $m \in \Nset$ denote by 
$R_q[M_n]^{\geq m}$ the sum of the components of 
$R_q[M_n]$ of degree $\geq m$. We call an automorphism 
$\Phi$ of $R_q[M_n]$ {\em{unipotent}} if 
\[
\Phi(x_{kl}) - x_{kl} \in R_q[M_n]^{\geq 2}. 
\]   
Launois and Lenagan proved \cite[\S 1.7]{LaL2} 
that \cjref{LL-cnj} follows from the following conjecture:

\bcj{un} For all base fields $\KK$, $q \in \KK^*$ 
not a root of unity, and $n \in \Zset_+$ every unipotent 
automorphism of $R_q[M_n]$ is equal to the identity 
automorphism.
\ecj
We prove \cjref{un} and thus obtain the validity of \cjref{LL-cnj}.
The key ingredient in the proof is the method of rigidity of quantum tori
which we developed in \cite{Y-ad}. The proof is in Section \ref{aut}.
Section \ref{qalg} gathers the needed 
facts for quantum matrices and quantum tori.
%\medskip
%\\
%\noindent
%{\bf Acknowledgements.} 
%The author was supported in part by NSF grant DMS-1001632.
%%%%%%%%%%%%%%%%%%%%%%%%%%%%%%%%%%%%%%%%%%%%%%%%%%%%%
\sectionnew{Quantum matrices and quantum tori}
\lb{qalg}
%%%%%%%%%%%%%%%%
\subsection{}
\label{2.1}
Throughout the paper, for all integers $k, l \in \Zset$, $k \leq l$,
we will denote $[k,l] := \{ k, \ldots, l \}$.

Let $N\in \Zset_+$. A multiplicatively skew-symmetric matrix
is a matrix $\qb=(q_{ij}) \in M_N(\KK^*)$ such that 
$q_{ii} =1$ for all $i \in [1, N]$ and $q_{ij} q_{ji} = 1$,
for all $1 \leq i < j \leq N$. To it one associates the 
{\em{quantum torus}} $\TT_\qb$ which is the $\KK$-algebra with 
generators $Y_i^{\pm 1 }$, $i \in [1,N]$ and relations
\[
Y_i Y_j = q_{ij} Y_j Y_i, \; \; 
\forall 1 \leq i < j \leq N, \quad
Y_i Y_i^{-1}=Y_i^{-1} Y_i = 1, \; \; \forall i \in [1, N].
\] 
The quantum torus $\TT_\qb$ is called {\em{saturated}} if 
for $u \in \TT_\qb$ and $k \in \Zset_+$, 
$u^k \in Z(\TT_\qb) \Rightarrow u \in Z(\TT_\qb)$.
Define the multiplicative kernel of $\qb$ by 
\[
\Ker (\qb) = \Big\{ (m_1, \ldots, m_N) \in \Zset^N
\mid \prod_{j=1}^N q_{ij}^{m_j} = 1, \; 
\forall i \in [1,N] \Big\}. 
\] 
The quantum torus $\TT_\qb$ is saturated if and only if
\[
(\Ker(\qb)/k) \cap \Zset^N = \Ker(\qb), \; \; 
\forall k \in \Zset_+,
\] 
see \cite[\S 3.1]{Y-ad}. In other words, for all $f \in \Zset^N$ and $k \in \Zset_+$
\[
kf \in \Ker(\qb) \Rightarrow f \in \Ker(\qb).
\]
It is clear that, if $q_{ij}$, $1 \leq i < j \leq N$ 
generate a torsion-free subgroup of $\KK^*$, 
then $\TT_\qb$ is saturated.
A vector $(d_1, \ldots, d_N) \in \Zset_+^N$
will be called a {\em{degree vector}}. It gives rise to a $\Zset$-grading
on $\TT_\qb$ by setting $\deg Y_i^{ \pm 1} = \pm d_i$. 
Denote by $\TT_\qb^r$ the subspace of $\TT_\qb$ of degree $r \in \Zset$ 
and consider the completion 
\[
\wh{\TT}_{\qb, \db} : = 
\{ u_r + u_{r+1} + \ldots \mid r \in \Zset, 
u_k \in \TT_\qb^k \; \mbox{for} \; k \geq r \}. 
\]
For $r \in \Zset$, let 
$\wh{\TT}_{\qb, \db}^{\geq r} : = 
\{ u_r + u_{r+1} + \ldots \mid r \in \Zset, 
u_k \in \TT_\qb^k \}$.
We call a continuous automorphism $\phi$ of $\wh{\TT}_{\qb, \db}$ 
{\em{unipotent}} if 
\[
\phi(Y_i) - Y_i \in \wh{\TT}_{\qb, \db}^{\geq d_i +1}, \; \; 
\forall i \in [1,N].
\] 
Under this condition $\phi(Y_i) = (1 + u_i) Y_i$
for some $u_i \in \wh{\TT}_{\qb,\db}^{\geq 1}$. Furthermore, 
$\phi(Y_i^{-1})$ is given by 
\[
\phi(Y_i^{-1})= 
( (1+u_i) Y_i )^{-1} =
\sum_{r=0}^\infty (-1)^r Y_i^{-1} u_i^r \in \wh{\TT}_{\qb, \db}^{\geq - d_i}
\]
and all values of $\phi$ are explicitly expressed in terms of $u_1, \ldots, u_N$. 
We will say that a unipotent automorphism $\phi$ of $\wh{\TT}_{\qb, \db}$ 
is {\em{bi-finite}} if 
\[
\phi(Y_i), \phi^{-1}(Y_i) \in \TT_\qb, \; \; \forall i \in [1,N].
\]  
We refer the reader to \cite[Section 3]{Y-ad} for details 
on the above notions. 

We will need the following result from \cite{Y-ad}, see Theorem 1.2 there.

\bth{1-bua} Fix an arbitrary base field $\KK$, a multiplicatively 
skew-symmetric matrix $\qb \in M_N(\KK^*)$ for which the quantum torus
$\TT_\qb$ is saturated, and a degree vector $\db \in \ZNp$. 
For every bi-finite unipotent automorphism
$\phi$ of the completed quantum torus $\wh{\TT}_{\qb, \db}$, 
there exists an $N$-tuple
\[
(u_1, u_2, \ldots, u_N) \; \; \mbox{of elements of } \; \; 
Z(\TT_\qb)^{\geq 1}
\] 
such that $\phi(Y_i) = (1+u_i) Y_i$
for all $i \in [1,N]$, where 
$Z(\TT_\qb)^{\geq 1} := Z(\TT_\qb) \cap \wh{\TT}_{\qb, \db}^{\geq 1}$.
\eth
%%%%%%%%%%%%%%%%
\subsection{}
\label{2.2}
Returning to the algebras of quantum matrices $R_q[M_n]$, 
fix a positive integer $n$ and denote 
for brevity
\begin{equation}
\label{Nn}
N:= n^2.
\end{equation}
Define $X_1, \ldots, X_N \in R_q[M_n]$ by
\[
X_{k + (l-1)n} : = x_{kl} \; \; \mbox{for} \; \; k, l \in [1,n].
\]
The form of the commutation relations for $x_{kl}$ produces 
an iterated Ore extension presentation for $R_q[M_n]$
\begin{equation}
\label{Ore}
R_q[M_n] = \KK[X_1][X_2; \sigma_2, \delta_2] \ldots 
[X_N; \sigma_N, \delta_N].
\end{equation}
Denote by $R_j$ the $j$-th algebra in the chain, i.e., 
the subalgebra of $R_q[M_n]$ generated by $X_1, \ldots, X_j$. 
For all $j \in [1, N]$, $\sig_j$ is an automorphism 
of $R_{j-1}$ such that   
\begin{equation}
\label{a}
\sig_j(X_i) = q^{a_{ji}} X_i, \; \; \forall i < j  
\end{equation}
for some integers $a_{ji}$, the explicit form of which 
will not play any role. Moreover, 
$\delta_j$ is a 
locally nilpotent $\sigma_j$-derivation of $R_{j-1}$
and $\sigma_j \delta_j = q^{-2} \delta_j \sigma_j, \; 
\forall \, 2 \leq j \leq N$. Cauchon's procedure of deleting derivations \cite{Ca1} 
constructs a sequence of $N$-tuples $(X^{(k)}_1, \ldots, X^{(k)}_N)$
of the division ring of fractions $\Fract(R_q[M_n])$ for
$k= N+1, \ldots, 2$. The first one 
is given by 
\[
(X^{(N+1)}_1, \ldots, X^{(N+1)}_N) := 
(X_1, \ldots, X_N)
\]
and the others are constructed inductively by 
\begin{equation}
\label{new-x}
X^{(k)}_j := 
\begin{cases}
X^{(k+1)}_j, 
& \mbox{if} \; \; j \geq k 
\\
\sum_{r=0}^\infty \frac{(1- q^{-2})^{-r}}{[r]_{q^{-2}}!} 
\Big[ \delta_k^r \sig^{-r}_k \left(x^{(k+1)}_j \right) \Big]
\left(x^{(j+1)}_j \right)^{-r}, 
& \mbox{if} \; \; j < k
\end{cases}
\end{equation}
for $k =N, \ldots, 2$ (in terms of the standard notation for $q$-integers 
and factorials $[0]_q=1$, $[r]_q = (1-q^r)/(1-q)$ for $r>0$, and
$[r]_q! = [0]_q \ldots [r]_q$ for $r \geq 0$).
Denote the final $N$-tuple
\begin{equation}
\label{final}
\ol{X}_1:= X_1^{(2)}, \ldots, \ol{X}_N := X_N^{(2)} \in \Fract(R_q[M_n])
\end{equation}
and the subalgebra
\[
\TT = \lcor \ol{X}^{\, \pm 1}_1, \ldots, \ol{X}^{\, \pm 1}_N \rcor 
\subset \Fract(R_q[M_n]).
\]
Define the multiplicatively skew-symmetric matrix
\begin{equation}
\label{matrix}
\qb = (q_{ij}) \in M_N(\KK^*) \; \;
\mbox{such that} \; \;  q_{ji} = q^{a_{ji}},
\; \; \forall \; 1 \leq i < j \leq N
\end{equation} 
for the integers $a_{ji}$ defined in \eqref{a}.
Cauchon proved in \cite{Ca1} the $\KK$-algebra isomorphism
\begin{equation}
\label{torus}
\TT_\qb \cong \TT, \; \; 
\mbox{where} \; \; Y_j \mt \ol{X}_j, \; 
1 \leq j \leq N.
\end{equation}
  
Define the map 
\begin{equation}
\label{mu}
\mu \colon [1,N] \to [1,2n-1] \; \; 
\mbox{by} \; \; 
\mu( k + (l-1)n) := n + l -k, \; \; 
\forall k,l \in [1,n].
\end{equation}
The meaning of this map is as follows. The algebra of quantum matrices 
$R_q[M_n]$ is isomorphic to the quantum Schubert cell 
subalgebras $\UU_q^\pm[w]$ of $\UU_q({\mathfrak{sl}}_{2n})$ for the 
following choice of Weyl group element 
\begin{equation}
\label{w}
w = (s_n \ldots s_1) \ldots (s_{2n-1} \ldots s_n) \in S_{2n} 
\end{equation}
which is in reduced form. Here $s_1, \ldots, s_{2n-1}$
denote the simple reflections of the symmetric group $S_{2n}$. We refer the 
reader to \cite[\S 2.1]{CM} and \cite[Section 4]{Y-gl} for details.
Within this framework the generators $X_1, \ldots, X_N$ of $R_q[M_n]$ 
are matched with the simple reflections in the expression 
\eqref{w} read from left to right.
The map $\mu$ simply encodes 
this matching. Define the following successor function 
$s \colon [1,N] \to [1,N] \cup \{ \infty \}$ attached to $\mu$:
\[
s(i) = \min \{ j > i \mid \mbox{such that} \; \; 
\mu(j) = \mu(i) \}
\; \; \mbox{if such a $j$ exits, otherwise} \; \; 
s(i) = \infty.
\]

Given two subsets $I = \{ k_1 < \ldots < k_r \}$ and $J = \{ l_1 < \ldots < l_r \}$
of $[1,n]$ having the same cardinality, define the quantum minor
\[
[I|J] = \sum_{w \in S_r} (-q)^{\ell(w)}
x_{k_1 l_{w(1)}} \ldots x_{k_r l_{w(r)}} \in R_q[M_n].
\]
We will use the following special quantum minors $\Delta_1, \ldots, \Delta_N \in R_q[M_n]$
(which combinatorially are the ones that contain consecutive
rows and columns, and  
either touch the right side or the bottom of an $n \times n$ 
matrix). Given $i \in [1,N]$, we represent it as $i = k+ (l-1)n$ 
for some $k,l \in [1,n]$ and define
\[
\Delta_i =
\begin{cases}
[ \{k, \ldots, n \} | \{l, \ldots, n+l-k \} ], & \mbox{if} \; \;  k \geq l
\\
[ \{k, \ldots, n- l +k \} | \{l, \ldots, n \} ], & \mbox{if} \; \; 
k < l. 
\end{cases}
\] 
The following theorem is due to Cauchon \cite[Proposition 5.2.1]{Ca2}. The above setting 
and the statement of this result follows the 
framework of \cite[Theorem 1.3]{GY} which establishes such a fact for all quantum 
Schubert cell algebras. The theorem below also follows from it and 
\cite[Lemma 4.3]{Y-gl}.

\bth{Ca1} For all $i \in [1,N]$ (recall \eqref{Nn}) the final $N$-tuple 
$(\ol{X}_1, \ldots, \ol{X}_N)$ of Cauchon's elements of $\Fract( R_q[M_n])$ 
is given by 
\begin{equation}
\label{de}
\ol{X}_i = 
\begin{cases} 
\Delta_i \Delta_{s(i)}^{-1}, & \mbox{if} 
\; \; s(i) \neq \infty
\\
\Delta_i, & \mbox{if} 
\; \; s(i) = \infty. 
\end{cases}
\end{equation}
\eth
%%%%%%%%%%%%%%%%%%%%%%%%%%%%%%%%%%%%%%%%%%%%
\sectionnew{Automorphisms of square quantum matrices}
\label{aut}
%%%%%%%%
\subsection{}
\label{3.1}
In this section we prove \cjref{un} which implies the
validity of the Launois--Lenagan conjecture. 
The proof of the following theorem is given in \S \ref{3.5}.

\bth{1} Let $\KK$ be an arbitrary base field, $q \in \KK^*$ 
an element which is not a root of unity, and $n \in \Zset_+$.
Every unipotent automorphism of $R_q[M_n]$ equals the 
identity automorphism.
\eth 

For $n \in \Zset_+$, we have the group embedding 
\[
\eta \colon \HH \rtimes \Zset_2 \hra \Aut R_q[M_n], \; \;
\eta(h, \ol{k}) (u) = h \cdot (\tau^k (u)), \; \; h \in \HH, k = 0,1
\] 
in terms of the notation from the introduction. Launois and Lenagan 
proved in \cite[Proposition 1.9]{LaL2} that every automorphism 
of $R_q[M_n]$ is a composition of an automorphism in the image of
$\eta$ and a unipotent automorphism. Hence \thref{1} implies 
the validity of \cjref{LL-cnj}.
\bth{LL-t} For all base fields $\KK$, $q \in \KK^*$ not a 
root of unity, and integers $n > 1$, the map 
$\eta \colon \KK \rtimes \Zset_2 \to \Aut R_q[M_n]$ is a group 
isomorphism. 
\eth

For completeness, we note that $R_q[M_1] = \KK [x_{11}]$ and 
$\Aut R_q[M_n]  \cong \KK \rtimes \KK^*$.  
%%%%%%%%%%%%%%%%%%%%%
\subsection{}
\label{3.1a}
Recall the definition \eqref{torus} of the quantum torus $\TT$ and 
the notation $N:=n^2$. 
Eq. \eqref{new-x} implies that the $\Nset$-grading of $R_q[M_n]$ from \S \ref{2.1}
extends to a $\Zset$-grading of $\TT$ by defining $\deg \ol{X}_i^{\, \pm 1} 
= \pm 1$ for all $i \in [1,N]$. Furthermore, for all intermediate steps of the 
Cauchon procedure
\[
\deg X^{(3)}_i = \ldots = \deg X^{(N)}_i = \deg X_i= 1, \; 
\; \forall i \in [1,N].
\]
Denote the corresponding graded subspaces of $\TT$ by $\TT^r$, $r \in \Zset$.
\thref{Ca1} implies that $\{\De_1, \ldots, \De_N\}$ is another set 
of independent generators of $\TT$. For this set of generators 
of $\TT$ consider the degree vector
$\db = (d_1, \ldots, d_N)$ where $d_i$ is the size of the quantum minor 
$\Delta_i$. If $i = k + (l-1) n \in [1,N]$ for $k, l \in [1,n]$, then
\begin{equation}
\label{d-vec}
d_i = n+ 1 - \max \{k,l\}.
\end{equation} 
By \thref{Ca1} the above $\Zset$-grading of $\TT$ coincides with the 
$\Zset$-grading from \S \ref{2.1} associated to this degree vector $\db$.
Denote by $\wh{\TT}$ the corresponding completion
as defined in \S \ref{2.1}. Let 
$\wh{\TT}^{\geq r} := \{ u_r + u_{r+1} + \ldots \mid u_k \in \TT^k \; \; \mbox{for} \; \; 
k \geq r \}$, $r \in \Zset$. For every subalgebra $R$ of $\wh{\TT}$
and $r \in \Zset$ denote $R^{\geq r} := R \cap \wh{\TT}^{\geq r}$ and 
$R^r := R \cap \TT^r$. 

It follows from \thref{Ca1} that for all $i \in [1,N]$
\[
\De_i = \ol{X}_i \ldots \ol{X}_{s^{n(i)}(i)}, 
\]
where $n(i)$ is the largest natural number such that $s^{n(i)}(i) \neq \infty$.
In particular, the generators $\De_1, \ldots, \De_N$ of $\TT$ satisfy the relations 
\[
\De_i \De_j = q'_{ij} \De_j \De_i,
\] 
where 
\[
q'_{ij} := \prod_{k =0}^{n(i)} \prod_{l=0}^{n(j)} q_{s^k(i), s^l(j)}.
\]

Analogously to \cite[\S 4.2]{Y-ad} we define an injective map 
from the group of unipotent automorphisms of $R_q[M_n]$ to 
the set of bi-finite unipotent automorphisms of $\wh{\TT}$.
Here is the general form of this correspondence:

\bpr{emb} Assume that $\TT_\qb$ is a quantum torus 
with generators $Y_1, \ldots, Y_N$ as in \S \ref{2.1} 
and $\db \in \Zset_+^N$ a degree vector. Let 
$R$ be a connected $\Nset$-graded subalgebra of $\TT_\qb$ 
(equipped with the $\Zset$-grading associated to $\db$)
such that $Y_i \in R$, $\forall i \in [1,N]$. 
For each automorphisms $\Phi$ of $R$ satisfying 
\begin{equation}
\label{Phi-unip}
\Phi(u) - u \in R^{\geq r + 1}, \; \; 
\forall u \in R^r, r \in \Nset,
\end{equation}
there exists a unique bi-finite unipotent automorphism $\phi$ of 
$\wh{\TT}_{\qb, \db}$ such that
$\phi(Y_i) = \Phi(Y_i)$, $\forall i \in [1,N]$. The assignment
$\Phi \mt \phi$ defines an injective homomorphism from the group 
of automorphisms $\Phi$ of $R$ satisfying \eqref{Phi-unip} 
into the group of unipotent automorphisms of $\wh{\TT}_{\qb, \db}$.
The image of this homomorphism is contained in the set of bi-finite unipotent 
automorphisms of $\wh{\TT}_{\qb, \db}$.
\epr 
\begin{proof} Let $\qb= (q_{ij})$ and $\db = (d_1, \ldots, d_N)$. 
We have $\Phi(Y_i) = Y_i + u'_i$, 
where $u'_i \in R^{\geq d_i +1} \subset \TT_\qb^{\geq d_i + 1}$,
$i \in [1,N]$.
Write $u'_i = u_i Y_i$ for some $u_i \in \TT_{\qb}^{\geq 1}$. Then 
$\Phi(Y_i) = (1+u_i) Y_i$ and 
$(1+u_i) Y_i (1 + u_j) Y_j = q_{ij} (1+u_j) Y_j (1+u_i) Y_i$, 
$\forall i,j \in [1,N]$, because $\Phi$ is an automorphism of $R$.
By \cite[Lemma 3.4]{Y-ad} there exists a unique unipotent 
automorphism $\phi$ of $\wh{\TT}_{\qb, \db}$ given by 
\begin{equation}
\label{ph}
\phi(Y_i) := (1+u_i) Y_i = \Phi(Y_i).
\end{equation}
Since $\TT_\qb$ is a domain, $\phi|_R = \Phi$.
Denote by $\psi$ the unipotent automorphism 
of $\wh{\TT}_{\qb, \db}$ associated by this construction to $\Phi^{-1}$. 
We have $(\psi \phi)|_R = (\phi \psi)|_R =\id$.
Hence $\psi = \phi^{-1}$ because $Y_1, \ldots, Y_N \in R$. 
This and Eq. \eqref{ph} imply that $\phi$ is a 
bi-finite unipotent automorphism of $\wh{\TT}_{\qb, \db}$.    
The map $\Phi \mt \phi$ is injective since $\Phi = \phi|_R$.

Assume that $\Phi_1$ and $\Phi_2$ are two automorphisms of $R$ satisfying
\eqref{Phi-unip}. Let $\phi_1$ and $\phi_2$ be the corresponding 
bi-finite unipotent automorphisms of $\wh{\TT}_{\qb, \db}$. 
Then $\phi_1 \phi_2(Y_i) = \phi_1 \Phi_2(Y_i) = 
\Phi_1 \Phi_2 (Y_i) \in R^{\geq d_i +1} \subset \TT_\qb^{\geq d_i +1},$ 
$\forall i \in [1,N]$ because ${\phi_1}|_R = \Phi_1$ and 
$\Phi_2(Y_i) \in R$. Therefore $\phi_1 \phi_2$ is also a bi-finite unipotent 
automorphism of $\wh{\TT}_{\qb, \db}$ and is precisely the automorphism 
associated to $\Phi_1 \Phi_2$.
\end{proof}

We apply \prref{emb} to the quantum torus $\TT$ with generating set 
$\{ \De_1, \ldots, \De_N\}$ and the degree vector $\db$ from Eq. \eqref{d-vec}.
The algebra $R_q[M_n]$ satisfies the conditions of the lemma since 
$\De_i \in R_q[M_n]$ and the $\Nset$-grading of $R_q[M_n]$ is the 
restriction of the $\Zset$-grading of $\TT$. Each unipotent automorphism 
$\Phi$ of $R_q[M_n]$ satisfies
\[
\Phi(u) - u \in R_q[M_n]^{\geq r+1}, \; \; 
\forall u \in R_q[M_n]^r, r \in \Nset
\]
because $R_q[M_n]$ is generated by $X_1, \ldots, X_N$. 
Thus, to each unipotent automorphism $\Phi$ of $R_q[M_n]$
\prref{emb} associates a bi-finite unipotent automorphism 
$\phi$ of $\wh{\TT}$ such that $\phi(\De_i) = \Phi(\De_i)$ 
for all $i \in [1,n]$. 

\bpr{1} Let $\Phi$ be a unipotent automorphism of $R_q[M_n]$, 
where $q \in \KK^*$ is not a root of unity.
Then there exist 
$u_1, \ldots, u_N \in Z(\TT)^{\geq 1}$ such that 
\[
\Phi(\De_i) = (1 + u_i) \De_i, \; \; \forall i \in [1,N]. 
\] 
\epr
\begin{proof}
Recall from \cite[\S 3.1]{Y-ad} that the property of a quantum torus being 
saturated does not depend on the choice of set of generators. 
The entries $q_{ij}$ of the the multiplicatively skew-symmetric 
matrix $\qb$ given by \eqref{matrix} are powers of $q$. 
They generate a torsion-free 
subgroup of $\KK^*$ because $q$ is not a root of unity. 
Thus $\TT_\qb$ is saturated, and due to the 
isomorphism \eqref{torus}, $\TT$ is also saturated.

Consider the bi-finite unipotent automorphism $\phi$ of $\wh{\TT}$ 
associated to $\Phi$ by the above construction. 
Applying \thref{1-bua}, we obtain 
\[
\phi (\De_i) = (1+ u_i) \De_i \; \; 
\mbox{for some} \; \; u_i \in Z(\TT)^{\geq 1}, \; \; 
\forall i \in [1,N].
\] 
Taking into account that $\Phi(\De_i) = \phi(\De_i)$ implies the 
statement of the proposition.
\end{proof}
%%%%%%%%%%%%%%%%%%%%
\subsection{}
\label{3.5}
Recall the definition of the surjective map 
$\mu \colon [1,N] \to [1, 2n-1]$ from Eq. \eqref{mu}.
The height one $\HH$-prime ideals of $R_q[M_n]$ 
are encoded in it as follows.  
For $i \in [1, 2n-1]$ denote   
\[
f(i) :=\min \mu^{-1}(i).
\]
In other words $f(i)$ records the first element 
of the preimage $\mu^{-1}(i)$. 
Launois, Lenagan, and Rigal 
proved \cite[Proposition 4.2]{LLR} that, in this 
framework, the height one $\HH$-prime ideals of 
$R_q[M_n]$ are 
\begin{equation}
\label{h1prim}
R_q[M_n] \De_{f(i)} \; \; \mbox{for} \; \; 
i \in [1, 2n-1].
\end{equation}
(In particular, all $\De_{f(i)}$ are normal elements 
of $R_q[M_n]$.) This is again a general property 
of all quantum Schubert cell algebras: all of their
height one torus-invariant prime ideals
are generated by a quantum minor 
which is a normal element \cite[Proposition 6.9]{Y-ams}.
Launois and Lenagan proved \cite[Proposition 1.10]{LaL2}
that every unipotent automorphism 
$\Phi$ of $R_q[M_n]$ fixes the height one prime ideals \eqref{h1prim}, i.e., 
fixes the elements $\De_{f(i)}$,
\begin{equation}
\label{fix}
\Phi(\De_{f(i)}) = \De_{f(i)}, \; \; 
\forall i \in [1,2n-1].
\end{equation}
\noindent
{\em{Proof of \thref{1}}}. The idea of the proof is similar
to the proof of \cite[Proposition 3.14]{Y-ad}.

Let $\Phi$ be a 
unipotent automorphism of $R_q[M_n]$ and $\phi$ 
the corresponding bi-finite unipotent automorphism 
of $\wh{\TT}$. \prref{emb} implies that 
\[
\Phi(\De_i) = (1 + u_i) \De_i \; \; 
\mbox{for some} \; \; u_i \in Z(\TT)^{\geq 1}, \; \; 
\forall i \in [1,N].
\]
The center of $\TT$ was described by 
Launois and Lenagan in \cite[Theorem 3.4]{LaL1}:
\[
Z(\TT) = \KK [ (\De_{f(1)} \De_{f(2n-1)}^{-1})^{\pm 1}, \ldots, 
(\De_{f(n-1)} \De_{f(n+1)}^{-1})^{\pm 1}, 
\De_{f(n)}^{\pm 1}].
\] 
Eq. \eqref{fix} and the fact that $\phi(\De_i) = \Phi(\De_i)$, 
$\forall i \in [1,N]$ imply 
\begin{equation}
\label{fix2}
\phi(u) = u, \forall u \in Z(\TT).
\end{equation}
Denote by $\psi$ the bi-finite unipotent automorphism 
of $\wh{\TT}$ corresponding to $\Phi^{-1}$. By Propositions \ref{pemb} and \ref{p1}, 
$\psi = \phi^{-1}$ and 
\[
\psi(\De_i) =
\Phi^{-1}(\De_i) = (1 + v_i) \De_i \; \; 
\mbox{for some} \; \; v_i \in Z(\TT)^{\geq 1}, 
\; \; \forall i \in [1,N].
\]
Applying Eq. \eqref{fix2} we obtain
\[
\De_i = \phi \psi (\De_i) = \phi( (1+v_i) \De_i) =
\phi(1+v_i) \phi(\De_i)= (1+v_i)(1+u_i) \De_i, 
\; \; \forall i \in [1,N].
\]
Because $\TT$ is a domain, we have 
\begin{equation}
\label{left-right}
(1+v_i) (1+u_i) = 1, \; \; \forall i \in [1,N].
\end{equation}
The condition $u_i, v_i \in Z(\TT)^{\geq 1}$ implies 
$u_i = v_i = 0$ for all $i \in [1,N]$. Indeed, 
if the highest terms of $1+ u_i$ and $1+ v_i$ are in degrees 
$r_i, t_i \in \Nset$, then the highest term of the left hand side 
of \eqref{left-right} is in degree $r_i + t_i$ because $\TT$ is a domain. 
The right hand side of \eqref{left-right} is in degree 0. 
Hence $r_i = t_i = 0$ and all elements $u_i$ and $v_i$ vanish. 
Therefore $\phi = \id$ and $\Phi= \id$. 
\qed 
%%%%%%%%%%%%%%%%%%%%%% References %%%%%%%%%%%%%%%%%%%%%%%%%%%%%%%%%%%%%%%

%%%%%%%%%%%%%%%%%%%%%%%%%%%%%%%%%%%%%%%%%%%%%%%%%%%%%%%%%%%%%%%%%%%%%%%%%%%%%%%
%%%%%%%%%%%%%%%%%%%%%%%%%%%%%%%%%%%%%%%%%%%%%%%%%%%%%%%%%%%%%%%%%%%%%%%%%%%%%%
\end{document}